\documentclass[12pt]{article}
\usepackage{amssymb}
\newcommand{\Rcc}{\mathbb{R}}
\newcommand{\Ccc}{\mathbb{C}}
\newcommand{\fr}{\frac}
\newcommand{\lb}{\label}
\newcommand{\ti}{\tilde}
\newcommand{\be}{\begin{equation}}
\newcommand{\ee}{\end{equation}}
\newcommand{\ba}{\begin{array}}
\newcommand{\ea}{\end{array}}
\newcommand{\beqa}{\begin{eqnarray}}
\newcommand{\beqay}{\begin{eqnarray*}}

\newcommand{\la}{\lambda}

\newcommand{\eeqa}{\end{eqnarray}}
\newcommand{\eeqay}{\end{eqnarray*}}
\newcommand{\ep}{\epsilon}
\newcommand{\Ac}{{\cal A}}

\newcommand{\Cc}{{\cal C}}
\newcommand{\Rc}{{\cal R}}
\newcommand{\Lc}{{\cal L}}
\newcommand{\Qc}{{\cal Q}}

\newcommand{\Ic}{{\cal I}}
\newcommand{\Kc}{{\cal K}}
\newcommand{\Omo}{\Omega}
\newcommand{\omo}{\omega}
\newcommand{\kd}{\delta}
\newcommand{\CD}{\Delta}
\newcommand{\epl}{\eta_+}
\newcommand{\emi}{\eta_-}
\newcommand{\ot}{\otimes}

\newcommand{\vae}{\varepsilon}
\begin{document}
\title{}
\author{}
\date{}

\noindent
{\Large \bf Non-Abelian Fractional Supersymmetry in Two Dimensions}

\vspace{1cm}
\noindent
{\footnotesize Haji AHMEDOV$^{a,}$}\footnote{E-mail address: 
hagi@gursey.gov.tr}
{\footnotesize and \"{O}mer F. DAYI$^{a,b,}$}\footnote{E-mail 
addresses: dayi@gursey.gov.tr and  dayi@itu.edu.tr.}

\vspace{10pt}

\noindent
$a)${\footnotesize \it Feza G\"{u}rsey Institute,}

\noindent
\hspace{3mm}{\footnotesize \it P.O.Box 6, 81220
\c{C}engelk\"{o}y--Istanbul, Turkey. }

\vspace{10pt}

\noindent
$b)${\footnotesize \it Physics Department, Faculty of Science and
Letters, Istanbul Technical University,}

\noindent
\hspace{3mm}{\footnotesize \it  80626 Maslak--Istanbul,
Turkey.}

\vspace{2cm}
\noindent
{\footnotesize {\bf Abstract.}} Non-Abelian fractional supersymmetry 
algebra in
two dimensions is introduced utilizing $U_q(sl(2,\Rcc ))$ at roots of
unity.
Its representations and the matrix elements are obtained.
The dual of it is constructed and the corepresentations are studied.
Moreover, a differential realization of the non-Abelian fractional
supersymmetry generators is given in the generalized superspace defined
by two commuting and two generalized Grassmann variables.
An invariant action under the fractional supersymmetry transformations
is given.

\newpage

\noindent 
{\bf 1. Introduction}

\noindent 
Supersymmetry is formulated in terms of $Z_2$ graded algebras
whose realizations can be obtained  by  Grassmann coordinates.
A step beyond supersymmetry is to consider $Z_3$ graded 
algebras\cite{ker}--\cite{rts}  
or the spaces given by  generalized Grassmann variables $\eta^p=0,$ 
where $p$ is a positive integer
and fractional supersymmetry generators 
\cite{abl}--\cite{ah},
which are defined 
as the $p^{\rm th}$ root of space--time translation operators.

Although, the latter approach can be studied
in terms of the
deformed algebras and  spaces  where the deformation parameter $q$
is a primitive root of unity, its 
group theoretical understanding was lacking. 
Recently, we proposed\cite{ho2} 
to study the two dimensional 
fractional supersymmetry
in terms of
a formulation of
the quantum Poincar\`{e}
group at roots of unity\cite{haci}.
This gave us the possibility of utilizing the well developed
representation theory techniques to study
the fractional supersymmetry in two dimensions.

In \cite{ho2}  two fractional supersymmetry generators which
are mutually commuting were considered. 
Here we generalize it to the fractional supersymmetry 
generators which are mutually noncommuting
utilizing the formulation given in \cite{ho1}.
Hence, we call it the 
non-Abelian fractional supersymmetry.
Our formulation straightforwardly leads to a 
differential realization of the fractional supersymmetry generators 
which would be very difficult, if not impossible, 
to find in terms of $q$--calculus.
The general theory of representations of $U_q(sl(2)$ at roots of unity
is well known\cite{ck}-\cite{bl} (and the references therein),
but for physical applications their explicit forms are needed as it is
discussed in the last section.

Another advantage of 
our formalism shows up when one deals with higher dimensions.
Although it is not presented here, 
higher dimensional fractional supersymmetry can be introduced
in terms of 
$SL_q(n, \Rcc )$ at roots of unity and its
subgroups by generalizing  our procedure of using 
$SL_q(2, \Rcc )$
to obtain the properties of the two dimensional
fractional supersymmetry, 
in a straightforward fashion.

After presenting the two dimensional
non-Abelian fractional supersymmetry
algebra $U_{FS}$
and its dual ${\cal A}_{FS},$ we deal with their
representations and corepresentations.
We study the $*$-representations
and in terms of them find differential
realizations of the non-Abelian 
fractional
supersymmetry generators in the superspace with
two commuting and 
two generalized Grassmann coordinates.
Then, we discuss how one can utilize these representations
in possible physical applications.
  
\vspace{1.5cm}

\noindent
{\bf 2. Non-Abelian fractional supersymmetry algebra and its dual  }

\vspace{.5cm}

\noindent
The two dimensional non-Abelian fractional supersymmetry algebra
 denoted $U_{FS},$ is the $*$--algebra generated by
$P_\pm,\ H,\ E_\pm ,\ K,$
 ($q^p=1,$ $p$ being an odd positive integer)
\be
\lb{f1}
{[ P_+ ,P_- ]} =   0,\ 
{[ P_\pm , H]}   =  \pm i P_\pm ,
\ee
\be
KE_\pm K^{-1} =q^{\pm 1}E_\pm ,\ \  \
[E_+,E_-]=\fr{K^2-K^{-2}}{q-q^{-1}},
\ee
\be
[K ,H]=0,\   [E_\pm ,H]=\pm \fr{i}{p} E_\pm ,
\ee
\be
\lb{f2}
P_\pm^*=P_\pm  ,\  H^*=H,\ E_\pm^*=E_\pm ,\  K^*=K ,
\ee
and $E_\pm$ are the $p^{\rm th}$ root of the
space-time translations $P_\pm :$  
\be
\lb{f3}
E_\pm^p=P_\pm. 
\ee
Moreover, we put the condition
\begin{equation}
\lb{f4}
K^p=1_U,
\end{equation}
where  $1_U$ indicates the unit element
of the algebra.

The basis elements of $U_{FS}$  are 
\be
\phi^{nmkrsl}\equiv E_-^nE_+^m K^k P_+^rP_-^sH^l ,
\ee
where $n,m,k=0,1,\cdots ,p-1,$ and $r,s,l$ are positive integers.

We can equip  $U_{FS}$ with the Hopf algebra structure 
\be
\begin{array}{lll}
\Delta (P_{\pm })=P_{\pm }\otimes 1_U+1_U\otimes P_{\pm },&
\varepsilon (P_{\pm })=0,&
S(P_{\pm})=-P_{\pm },\\
\CD (H) =H \ot 1_U   +1_U   \ot H, & \vae (H)=0, & S(H)=-H, \\
\Delta (E_\pm )=E_\pm \otimes K+K^{-1}\otimes E_\pm , &
\varepsilon (E_\pm )=0 , &
 S(E_\pm )=-q^{\pm 1}E_\pm , \\
\Delta (K)= K\otimes K, &
\varepsilon (K^{\pm 1})=1, &
S(K^{\pm 1})=K^{\mp 1}.
\end{array}
\ee

Let us present  the dual
of the non-Abelian fractional supersymmetry algebra $U_{FS}:$
It is the $*$--algebra
$\Ac_{FS}$ with $q^p=1$ and
generated by $z_\pm ,\ \la,\ \eta_\pm,\ \kd $ satisfying
\beqa
\eta _{-}\eta _{+}=q^2\eta _{+}\eta _{-}, & \eta _{\pm }\delta 
=q^2\delta \eta _{\pm }, &  \\
\eta_\pm^p=0 ,& \kd^p =1_A , &  \\
\eta_\pm^*=\eta_\pm,  & \kd^*=\kd  ,& \\
z_\pm^*=z_\pm, & \la^*=\la ,&  
\eeqa
where $z_\pm ,\ \la $ commute with the others and
 $1_A$ is the unit element of $\Ac_{FS} .$
Its Hopf algebra structure is given by the coproducts 
\be
\lb{ftra}
\begin{array}{lcl}
\CD \kd & = & \kd \otimes \kd + q^{-2}\kd^{-1}\eta_+^2\otimes
\eta_-^2\kd +(1_A+q^{-2})\eta_+\otimes \eta_- \kd, \lb{cp1} \\
\CD \eta_+& = & \epl \ot 1_A +\kd \ot \epl + (1_A+q^2)\epl \ot \epl \emi
+q^{-2}\kd^{-1}\epl^2 \ot (1_A+q^2 \epl \emi )\emi, \lb{cp2}  \\
\CD \eta_-& = & \emi \ot 1_A +\kd^{-1} \ot \emi +
\sum_{k=1}^{p-2} 
(-1)^k q^{-k(k+1)} \kd^{-k-1} \epl^k \ot \emi^{k+1},   \\
\CD \la  & = & \la \ot 1_A +1_A \ot \la ,\\
\CD z_+ & = & z_+ \ot 1_A +e^\la  \ot z_+ +
\sum_{k=1}^{p-1}\fr{q^{k^2}}{[k]![p-k]!}
\epl^{p-k} \kd^{k}e^{n\la /p} \ot (-q^2 \epl\emi ; q^{2})_{(p-k)}
\epl^k, \\
\CD z_- & = &  z_- \ot 1_A +e^{-\la } \ot z_- +
\sum_{k=1}^{p-1}\fr{q^{-k^2}}{[k]![p-k]!}
\emi^{p-k} \kd^{-k}e^{-n\la /p} (-\epl\emi ; q^{-2})_k \ot \emi^k,  
\end{array}
\ee
the antipodes 
\[
\begin{array}{ll}
S(\kd )  =  \kd^{-1} (1_A+q^{-2} \epl \emi )(1_A+\epl \emi ), &
S(\eta_\pm )  =  -\kd^{\mp 1} \eta_\pm , \\
S (\la )=-\la , &
S(z_\pm )   =  - z_\pm , \\
\end{array}
\]
and the counits
\[
\begin{array}{llll}
\ep (\kd )  =  1, & \ep (\eta_\pm )  =  0, &
\ep (\la ) =-\la , &  \ep (z_\pm )=0.
\end{array}
\]
We use the notation
\[
(a;q)_k \equiv \prod_{j=1}^k (1-aq^{j-1}) ,
\]
the symmetric $q$--number
\[
[n] =\fr{q^n-q^{-n}}{q-q^{-1}}
\]
and the $q$--factorial $[n]!=[n][n-1]\cdots [1].$

Any element of $\Ac_{FS}$ can be written as
\[
\sum_{n,m,k=0}^{p-1} f_{nmk}(z_+,z_-,\la )\epl^n \emi^n\la^k \in \Ac_{FS},
\]
where $ f_{nmk}(z_+,z_-,\la )$ are
infinitely differentiable functions on $\Rcc^3$.
Definition of the Hopf algebra operations
on these functions were given in \cite{ho1}.
Hence, a local basis of $\Ac_{FS}$ can be given by
\be
a^{nmktsl}\equiv \eta_-^n \eta_+^m \zeta (k,\kd ) z_+^t z_-^s \la^l ,
\ee
where $n,m,k=0,1,\cdots ,p-1;\ t,s,l$ are positive integers and
we defined
$$
\zeta (m, \kd )
\equiv \frac 1p\sum_{n=0}^{p-1}q^{-nm}\delta ^n.
$$
Observe that
\beqa
\langle \phi^{nmktsl},
a^{n^\prime m^\prime k^\prime t^\prime s^\prime l^\prime } \rangle
&  = &  i^{n+m+t+s+l} q^{\frac{n-m}{2}-nm}t!s!l![n]![m]! \nonumber \\
& &  \delta_{nn^\prime } \delta_{mm^\prime }
\delta_{tt^\prime }\delta_{ss^\prime } \delta_{ll^\prime }
\delta_{k+n+m, k^\prime }, \lb{dr}
\eeqa
are the duality relations between $\Ac_{FS}$ and $U_{FS}.$

\vspace{1cm}

\noindent
{\bf 3. Pseudo-unitary, irreducible corepresentations of  $\Ac_{FS}$}

\vspace{.5cm}

\noindent
Let $C_0^\infty (\Rcc )$ be the space of
all infinitely differential functions  of
finite support in $\Rcc$
and $P(t)$ denote the algebra of polynomials in $t$
subject to the conditions $t^p=1$ and $t^*=t,$
i.e. any element of $P(t)$ can be written as
$a(t)\equiv \sum_{n=0}^{p-1}a_nt^n.$
The irreducible representation of $U_{FS}$
in $C_0^\infty (\Rcc ) \times P(t) $ is defined by
the linear map
$$
\pi_{\la_\pm}(U_{FS} ):  C_0^\infty (\Rcc ) \times P(t)\rightarrow
C_0^\infty (\Rcc ) \times P(t)
$$
given as 
\beqa
\pi_{\la_\pm} (E_+ ) f(x) a(t) & = &
\la_+^{1/p} e^{x/p} t f(x)a(t) ,\nonumber \\
\pi_{\la_\pm} (E_- ) f(x) a(t) & = &
 e^{-x/p}  f(x)\sum_{n=0}^{p-1}M_n a_nt^{n-1}
 ,\nonumber \\
\pi_{\la_\pm} (P_\pm ) f(x) a(t) & = &
\la_\pm e^{\pm x} f(x)a(t),  \lb{cc} \\
\pi_{\la_\pm} (H ) f(x) a(t) & = &
-i \fr{d}{dx} f(x)a(t), \nonumber \\
\pi_{\la_\pm} (K ) f(x) a(t) & = &
 f(x)a(qt), \nonumber   
\eeqa
where we used the notation
\[
M_0=\la_+^{-1/p} ,\  M_n=\la_+^{-1/p}
\{ \la_+ -[n][n-1]\} \ {\rm for}\ n\neq 0,
\]
and
\[
\la_-=\prod_{n=0}^{p-1}M_n .
\]

Let us introduce the following hermitian forms
for the space $C_0^\infty (\Rcc ) \times P(t),$
\beqa
(f_1,f_2) & = & \int_{-\infty}^{+\infty}dxf_1(x)
\overline {f_2(x)} ,\lb{h1} \\
(a_1,a_2) & = & \Phi \left( a_1(t)a^*_2(t) \right), \lb{h2}
\eeqa
where
\be
\Phi(t^s)=\kd_{s,0{\rm (mod\ p)}}.
\ee

$C_0^\infty (\Rcc ) $ endowed with the norm induced by (\ref{h1})
leads to the Hilbert space of the square integrable
functions on $\Rcc .$ On the other hand $P(t)$ with the
norm $||a||^2 \equiv (a,a)$ is the pseudo-Euclidean space
with $\fr{p+1}{2}$ positive and $\fr{p-1}{2}$ negative
signatures\cite{ho1}. Now, one can verify that 
$\pi_{\la_\pm}$ defines pseudo-unitary,
irreducible $*$--representation of $U_{FS} $ for real $\la_\pm .$

The irreducible
corepresentation of $\Ac_{FS}$ can be
found in terms of
the duality relations (\ref{dr})
and the representation (\ref{cc}) of $U_{FS}$ as 
\be
\lb{tir}
T_{\la_\pm}(f(x)a(t)) = 
\sum_{n,m,k=0}^{p-1}\sum_{t,s,l=0}^{\infty}
\fr{a^{nmktsl} \pi_{\la_\pm}(\phi^{nmktsl})f(x)a(t) }
{\langle \phi^{nmktsl},a^{n m k t s l } \rangle},
\ee
which is pseudo-unitary for real $\la_\pm .$

Consider the Fourier transform of $f(x)\in C_0^\infty (\Rcc ) $
\be
F(\nu )  =  \int_{-\infty}^{+\infty} f(x) e^{\nu x} dx.
\ee
This integral converges for any complex $\nu .$
$F(\nu)$ is an analytic function and moreover, satisfies
\be
|F({\rm Re }\ \nu  +i{\rm Im }\ \nu) |< \omo e^{c|{\rm Re}\ \nu |},
\ee
for some real constants $\omo$ and $c.$
Then we can write the  inverse transform as
\be
f(x)  =  \fr{1}{2\pi i} \int_{c-i\infty }^{c+i\infty}
F(\nu )e^{-\nu x} d\nu .
\ee

The Fourier transform of $T_{\la_\pm}$ (\ref{tir})
yields the pseudo-unitary corepresentation 
in the space of functions $F(\nu )a(t)$ as
\be
\Qc_{\la_\pm}\left( F(\nu )t^k\right) = 
\int_{c-i\infty }^{c+i\infty} d\mu
\sum_{l=0}^{p-1} \Qc_{kl}^{\la_\pm}(\nu ,\mu , g)
F(\mu )t^l,
\ee
where $g\equiv (g_0;\ g_p)\equiv ( z_+,z_-,\la ;\ \epl ,\emi ,\kd).$
The kernel $\Qc_{kl}^{\la_\pm}$ is
\beqa
{\rm for}\ l\geq k, & & \lb{fa1} \\
\Qc_{kl}^{\la_\pm}(\nu ,\mu ,g )  & = &
 \Kc_{l-k}^{\la_\pm}(\nu ,\mu ,g_0) \kd^l\epl^{l-k}\Omo_{k,l}(\xi )
+ \Kc_{l-k-p}^{\la_\pm}(\nu ,\mu ,g_0) \emi^{p+k-l} \kd^l
\ti{\Omo}_{k,l} (\xi),    \nonumber \\
{\rm for}\ l< k, & &    \lb{fa2}  \\
Q_{kl}^{\la_\pm}(\nu ,\mu ,g )  & = &
 \Kc_{p+l-k}^{\la_\pm }(\nu ,\mu ,g_0) \kd^l\epl^{p+l-k}\Omo_{k,p+l}(\xi )
+ \Kc_{l-k}^{\la_\pm}(\nu ,\mu ,g_0) \emi^{k-l} \kd^l
\ti{\Omo}_{k,p+l} (\xi),   \nonumber
\eeqa
where we introduced,
in terms of
$\xi= q\epl \emi ,$ 
\be
\Omo_{k,l}(\xi )=
\sum_{m=0}^{p+l-k-1} 
\fr{i^{k-l}(-1)^{m+k+l}q^{(l-k)(l+1/2)-m(k+l)}
\la_+^{(m-k+l)/p}}{[m]![m-k+l]!}
\left( \prod_{s=0}^m M_{s+l} \right) \xi^m ,
\ee

\be
\ti{\Omo}_{k,l}(\xi )=
\sum_{m=0}^{k-l-1} 
\fr{i^{l-k-p}(-1)^{m+1}q^{l(k-l)+m(k+l)}
\la_+^{m/p}}{[m]![m+p+k-l]!}
\left( \prod_{s=1}^{p+k-l+m} M_{s-p+l} \right) \xi^m .
\ee
The functions $\Kc^{\la_\pm}_s$ are 
\be
\lb{sf}
\Kc_s^{\la_\pm}(\nu ,\mu ,g_0) = \fr{1}{2\pi  i}
e^{\mu \la_\pm } \int_{-\infty}^{+\infty}  e^{ ir
( e^xz_+ + e^{-x}z_-)+
x(\nu-\mu+ s/p)} dx.
\ee

By utilizing
the analog of  polar coordinates $\rho >0,\ \beta \in \Rcc $,
the pseudo-Euclidean plane
defined by the axis $z_-=0$ and $z_+=0$ can be studied
in terms of the quadrants
\beqa
{\rm Quad. 1}:&  z_+z_->0 ,\ z_\pm =\fr{1}{2}\rho e^{\pm \beta }, &
{\rm Quad. 2}:\ \ z_+z_-<0 ,\ z_\pm
=\pm \fr{ 1}{2}\rho e^{\pm \beta },  \nonumber   \\
{\rm Quad. 3}:& z_+z_->0,\ z_\pm =\fr{-1}{2}\rho e^{\pm \beta },  &
{\rm Quad. 4}:\ \  z_+z_-<0,\ z_\pm =\mp \fr{ 1}{2}\rho e^{\pm \beta }. 
\nonumber
\eeqa
In these quadrants (\ref{sf}) will lead to the
Hankel functions $H^{(1)}_\nu ,\  H^{(2)}_\nu$ or cylindrical
functions of imaginary argument $K_\nu :$ 
\beqa
{\rm Quad. 1:} &
\Kc_s^{\la_\pm}(\nu ,\mu ,g_0)   =
\fr{1}{2} e^{(\mu -\nu -s/p)( \beta + \fr{\pi i}{2} )
+\mu \la_\pm  } H^{(1)}_{\mu-\nu -s/p} (r\rho )  , & \nonumber \\
{\rm Quad. 2:} & \Kc_s^{\la_\pm}(\nu ,\mu ,g_0)   = 
\fr{1}{2} e^{(\mu -\nu -s/p)( \beta -\fr{\pi i}{2} )
+\mu \la_\pm   }H^{(2)}_{\mu-\nu -s/p} (r\rho ) , & \nonumber  \\
{\rm Quad. 3:} & \Kc_s^{\la_\pm}(\nu ,\mu ,g_0)  = 
\fr{1}{\pi i} e^{(\mu -\nu -s/p)( \beta +\fr{\pi i}{2} )
+\mu \la_\pm  } K_{\mu-\nu -s/p} (r\rho )  , &\nonumber \\
{\rm Quad. 4:} & \Kc_s^{\la_\pm}(\nu ,\mu ,g_0)   = 
\fr{1}{\pi i} e^{ (\mu -\nu -s/p)( \beta -\fr{\pi i}{2} )
+\mu \la_\pm  } K_{\mu-\nu -s/p} (r\rho )  , & \nonumber 
\eeqa
with the condition $-1< {\rm Re} (\nu -\mu + s/p)<1$\cite{vk}.

\vspace{1cm}

\noindent
{\bf 4. Quasi-regular corepresentation
of $\mathbf \Ac_{FS},$  $\mathbf *$--representation
of the non-Abelian  fractional supersymmetry algebra
and a differential realization}

\vspace{.5cm}

\noindent
The comultiplication
\be
\lb{de}
\CD :\Ac \rightarrow \Ac_{FS} \ot \Ac
\ee
defines the pseudo-unitary left quasi-regular corepresentation of
$\Ac_{FS}$ in its subspace $\Ac$ consisting of
the finite sums
$$
X=\sum_s a_s(\epl ,\emi )f_s (z_+,z_-)
$$
where
$a_s(\epl ,\emi )$ are polynomials in $\epl ,\emi$ and
$f_s (z_+,z_-) \in C_0^\infty (\Rcc^2) .$
The space $\Ac$  can be endowed with 
the hermitian form 
\be
\lb{222}
(X,Y)_E=\Ic_E(XY^*),
\ee
$X,Y\in\Ac$ and
the linear functional
$\Ic_E :\ \Ac \rightarrow \Ccc$
\be
\lb{333}
\Ic_E (X) =\sum_s \Ic (a_s) \Ic_C(f_s)
\ee
was shown to be
the left invariant integral\cite{ho1} in terms of the integrals
on generalized superspace\cite{rts},\cite{abl},\cite{dur}--\cite{ah},
\cite{bf}--\cite{is}
\beqa
\Ic (\epl^n \emi^m ) & = & q^{-1}\kd_{n, p-1} \kd_{m, p-1}, \\
\Ic_C(f_s)& = & \int_{-\infty}^{+\infty}dz_+dz_-f_s(z_+,z_-).
\eeqa

The right representation of the non-Abelian fractional
supersymmetry algebra $U_{FS}$
corresponding to the quasi-regular
representation (\ref{de}),
\be
\lb{rr}
\Rc (\phi )X= (\phi \ot id)\CD( X),
\ee
$\phi \in U_{FS},$ is a $*$--representation
\[
(\Rc (\phi )X,Y)_E =(X,\Rc (\phi^* )Y)_E ,
\]
due to the fact that
 the hermitian form (\ref{222}) is defined in terms of
the left
invariant integral (\ref{333}).

The right representations
on the variables $\eta_\pm$ and $f(z_+,z_-)$
can explicitly be written as
\be
\lb{reo}
\begin{array}{ll}
\Rc (E_+)\epl^n  =
iq^{1/2}[n]\epl^{n-1} + iq^{1/2-n}[2n]\emi \epl^{n},&
\Rc (E_+)\emi^n   =  - iq^{-1/2}[n]\emi^{n+1}, \\
\Rc (E_-)\epl^n  =    0, &
\Rc (E_-)\emi^n  =  iq^{-1/2} [n]\emi^{n-1},   \\
\Rc (K)\eta_\pm^n  =   q^{\pm n} \eta_\pm^n, &
\Rc (H)\eta_\pm^n=\pm \fr{in}{p}\eta_\pm^n,  \\
\Rc (P_\pm)\eta_\pm^n=0, &
\Rc (K)z_\pm =  z_\pm , \\
\Rc (E_\pm) f(z_+,z_-)  =  \fr{iq^{\pm 1/2}}{[p-1]!}\eta_\pm^{p-1}
\fr{df(z_+,z_-)}{dz_\pm }, 
&  \Rc (P_\pm)(z_+,z_-)=
i\fr{df(z_+,z_-)}{dz_\pm}, \\
\Rc (H) f= iz_+\fr{df}{dz_+}-iz_-\fr{df}{dz_-}. &
\end{array}
\ee

The  relations satisfied by
the right representation ${\cal R}$ 
\beqa
{\cal R}( \phi \phi^\prime ) & = & {\cal R}(\phi^\prime )
 {\cal R}(\phi ) , \nonumber\\
{\cal R}( E_\pm ) (X Y ) & = & {\cal R}( E_\pm )X {\cal R}( K ) Y
+ {\cal R}(K^{-1})X {\cal R}(E_\pm )Y,   \nonumber\\
{\cal R}(K ) (XY) & = & {\cal R}(K )X
{\cal R}(K )Y , \nonumber\\
\Rc (H) (XY) & = & \Rc (H)X Y+ X\Rc (H) Y , \nonumber \\
\Rc (P_\pm)(XY) & = & \Rc (P_\pm)X Y+ X\Rc (P_\pm) Y , \nonumber
\eeqa
permit us to define the action of an arbitrary
operator ${\cal R}(\phi )$ on any function in $\Ac.$

The quantum algebra which we deal with possesses two 
Casimir elements
\beqa
C_1 & = & E_-E_++\fr{(qK-q^{-1}K^{-1})^2}{(q^2-q^{-2})^2}, \lb{c1} \\
C_2 & = & P_+P_- . \lb{c2} 
\eeqa
As the complete set
of commuting operators we can choose $\Rc (C_1)$,  $\Rc (C_2)$
$\Rc (H),$  $\Rc (K ) $ and
$\Lc (H),$  $\Lc (K ) $ where
$\Lc (\phi )$ is the left representation of the element $\phi$
defined similar to (\ref{rr}) with the interchange of $\phi$
with the identity $id.$ Observe that
$\Lc (H) X=0$ and  $\Lc (K )X =X$ for any $X\in \Ac ,$
thus, in the space $\Ac$
the matrix elements can
be labeled as $D^{\la_\pm}_{n\nu ,m\mu }.$
Indeed, in terms of the kernel $\Qc^{\la_\pm}_{mn}$
 (\ref{fa1}) one observes that
\[
D_{n\nu ,00} = \Qc_{0n}^{\la_\pm}(\nu ,0 ,g ),
\]
$n \in [0,p-1],$ satisfy
\beqa
\Rc (K ) D_{n\nu ,00} & = & q^n D_{n\nu ,00}  \nonumber \\
\Rc (H ) D_{n\nu ,00} & = & -i(\nu+n/p) D_{n\nu ,00}  \nonumber \\
\Rc (C ) D_{n\nu ,00} & = & \la_+ D_{n\nu ,00}  \nonumber \\
\Rc (E_+)D_{n\nu ,00} & = & \la_+^{1/p}
D_{n + 1\nu ,00}, \nonumber \\
\Rc (E_-)D_{n\nu ,00} & = & M_n
D_{n- 1\nu ,00}, \nonumber \\
\Rc (P_\pm)D_{k\nu ,00} & = & \la_\pm
D_{k\nu \pm 1,00} ,\nonumber
\eeqa
where we introduced the notation
$ D_{p\nu ,00} \equiv  D_{0\nu +1 ,00} $
and
$ D_{-1\nu ,00} \equiv  D_{p-1\nu -1 ,00} .$

To write field theory actions in 
the generalized superspace given by $\eta_\pm,\ z_\pm ,$
one should  be equipped with ($q$-)differential realizations
of  the fractional supersymmetry generators, $E_\pm .$ When 
the functions of only  $\epl ,z_+$ or $\emi ,z_-$ 
are considered,
one can use either the algebraic or the group theoretical
properties to find the  differential
realizations of the generators \cite{abl},\cite{dur}-\cite{ho1}
in terms of the $q$--derivatives $D_\pm^q$ satisfying
\be
D^q_\pm \eta^n_\pm =\fr{1-q^n}{1-q} \eta_\pm^{n-1}. 
\ee

When we deal with 
the generalized superspace $\eta_\pm,\ z_\pm ,$
using the algebraic properties ($q$-calculus) to find the differential
realizations
of the non-Abelian fractional supersymmetry 
generators is hopeless. However,
the right representation obtained in (\ref{reo}) can be used to write
the operators corresponding to
 $E_\pm $ 
in the generalized superspace given by
$\eta_\pm,\ z_\pm$ with the ordering
\[
F(\eta_\pm,z_\pm)=\sum_{m,n}f_{mn}(z_+,z_-)\epl^n\emi^m ,
\]
in terms of the operators
\[
\hat{\eta}_\pm F=F\eta_\pm,
\]
and the $q$--derivatives $D_\pm^q$ satisfying
\be
\begin{array}{c}
D^q_\pm \hat{\eta}_\pm -q \hat{\eta}_\pm D^q_\pm =1, \\
D^q_\pm \hat{\eta}_\mp = q^{\pm 2} 
\hat{\eta}_\mp D^q_\pm .
\end{array}
\ee
This can be achieved by using the well known
mutually commuting ``dilatation" operators
$T_\pm$ and their inverses
\be
\begin{array}{c}
T_\pm =D_\pm^q\eta_\pm-\eta_\pm D_\pm^q=1 -(1-q)\hat{\eta}_\pm D_\pm^q,\\
 T^{-1}_\pm =1-(1-q^{-1})\hat{\eta}_\pm D_\pm^{q^{-1}},
\end{array}
\ee
satisfying
\be
\begin{array}{c}
T_\pm\eta_\pm^n=q^n\eta_\pm^n ,\\
T_\pm\eta_\mp^n=\eta_\mp^n .
\end{array}
\ee
Differential realizations of the fractional 
supersymmetry generators in the generalized superspace are
\beqa
\Rc (E_+ ) & = & \fr{iq^{1/2}}{1+q}T_-(D_+^{q^{-1}}+qD_+^q)
+\fr{iq^{1/2}}{q^2-q^{-2}}\hat{\eta}_-T^{-1}_-T^3_+
+ \fr{iq^{-3/2}}{q^2-q^{-2}}\hat{\eta}_-T^{-1}_-T^{-1}_+
 \nonumber \\
 & &  - \fr{iq^{-1/2}}{q-q^{-1}}\hat{\eta}_-T_-T^{-1}_+
+\frac{iq^{1/2}  }{[p-1]!}
\eta^{p-1}_+ \frac{d}{dz_+ },\lb{sps} \\
\Rc (E_- ) & = & \fr{iq^{-1/2}}{1+q}T_+^{-1}\left(
D_-^{q^{-1}}+qD_-^q \right)
+\frac{iq^{-1/2} }{[p-1]!}
\eta^{p-1}_- \frac{d}{dz_- } \lb{spr} .
\eeqa
In terms of the scalar product (\ref{333}) we can define
the involutions on the $q$--derivatives as
\beqa
{D_\pm^q}^* & = & -qD_\pm^q ,\\
{D_\pm^{q^{-1}}}^* & = &-q^{-1}D_\pm^{q^{-1}} ,
\eeqa
and verify that the realizations
(\ref{sps})-(\ref{spr}) satisfy
the involution conditions
\[
\Rc (E_\pm)^* =\Rc (E_\pm),
\]
as they should. Moreover, $\Rc (K)=T_-^{-1}T_+$ can easily be observed.

\vspace{1cm}

\noindent
{\bf 5. Discussions}

\vspace{.5cm}

\noindent
The formalism 
presented here can directly be used to define actions 
which are invariant under the transformations
of the generalized superspace,
(\ref{ftra}),
as far as we know the Casimir elements of the algebra.

Let ${\cal C}$ be any function of the Casimir elements $C_1,$
 $C_2,$ (\ref{c1})--(\ref{c2}), and the identity.
The action 
\be
\lb{ina}
S[\Phi]=\Ic_E (\Phi^*  \Rc ({\cal C})\Phi ),
\ee
where $\Phi (z_\pm,\eta_\pm )$  is
any function on
the generalized superspace, 
is invariant under the fractional
supersymmetry transformations (\ref{ftra}):
Being  the Casimir
element,  ${\cal C}$ satisfies
\[
\CD \Rc (\Cc )=\left( id \ot \Rc (\Cc )\right) \CD .
\]
Since  the Casimir operators
$\Rc (\Cc_{1,2})$  commute with 
the quasiregular representation $\CD$  of $\Ac_{FS},$
(\ref{de}), we get
\be
\begin{array}{lll}
S[\CD \Phi] & =
\left( id \ot \Ic_E \right) \left( \CD \Phi^* (id  \ot  \Rc (\Cc ) )
\CD \Phi \right)  & =  
\left( id \ot \Ic_E \right) \left( \CD \Phi^* \CD \left( \Rc (\Cc )
\Phi \right) \right)   \\
& =  \left( id \ot \Ic_E \right) \CD \left(  \Phi^*   \Rc (\Cc )
\Phi \right)  &= \Ic_E (\Phi^*  \Rc (\Cc ) \Phi )\\
& =S[\Phi] .&   
\ea
\ee
The choice of $\Cc$ 
depends on the physical system which we would like to deal with.
Equipped with these representations, 
(\ref{sps})--(\ref{spr}),
we hope that one can find a
fractional supersymmetric
action which can be useful to understand some physical problems like the
quantum Hall
effect where the $U_q(sl(2))$ at q roots of unity appears\cite{om}
in a natural way.

Fractional supersymmetric quantum mechanics\cite{fqm}
and field \mbox{theories\cite{abl}--\cite{azm},\cite{fik}}  were 
discussed in terms of q--calculus. Because of being a group
theoretical approach, we hope that our way of treating fractional
supersymmetry
will shed some light on how to  overcome the difficulties which show up in
defining fractional supersymmetric models.

\end{document}